\documentclass{amsart}
\usepackage{amssymb}
\usepackage{amsfonts}
\usepackage{latexsym}

\newtheorem{theorem}{Theorem}[section]
\newtheorem{lemma}[theorem]{Lemma}

\newtheorem{corollary}[theorem]{Corollary}
\theoremstyle{definition}
\newtheorem{definition}[theorem]{Definition}

\newtheorem{remark}[theorem]{Remark}

\newcommand{\Id}{\text{Id}}

\newcommand{\Ind}{\text{Ind}}

\newcommand{\eps}{\varepsilon}

\newcommand{\ben}{\begin{enumerate}}
\newcommand{\een}{\end{enumerate}}

\newcommand{\cO}{{\mathcal O}}
\newcommand{\cL}{{\mathcal L}}
\newcommand{\cF}{{\mathcal F}}

\newcommand{\BZ}{{\mathbb Z}}

\begin{document}

\title{A remark on cuspidal local systems}
\author{Victor Ostrik}
\address{Institute for Advanced Study, Einstein Drive, Princeton, NJ 08540}

\begin{abstract}
In this note we show that all reductive groups are
clean in characteristic $\ge 3$. In characteristic 2 there are
two cuspidal local systems (one for $F_4$ and one for $E_8$) which can not
be handled by our method.
\end{abstract}

\maketitle

\section{Introduction}
\subsection{}
Let $G$ be a connected reductive group over an algebraically closed
field of characteristic $p\ge 0$. In the study of the character sheaves on
the group $G$ (\cite{LuCS}) an important technical assumption is that
the group $G$ is clean. This means that all cuspidal local systems (see 
\cite{LuIC}) on all Levi subgroups of $G$ are clean (that is the IC-extension 
of these local systems coincides with the extension by zero). It is expected 
that any group $G$ is clean; equivalently any cuspidal local system is clean 
for any reductive group $G$. It is known \cite{LuCS} that the last assertion 
is equivalent to the similar assertion with $G$ assumed to be almost simple. 
It was shown in \cite{LuCS} that any cuspidal local system is clean if

1) $G$ is of classical type and $p$ is arbitrary;

2) $G$ is of type $E_6$ and $p=0$ or $p>2$;

3) $G$ is of type $G_2$, $F_4$, $E_7$ and $p=0$ or $p>3$;

4) $G$ is of type $E_8$ and $p=0$ or $p>5$.

In \cite{Sh1, Sh2} T.~Shoji used the Shintani descent theory to improve the
bounds for $p$ above for $G$ of types $G_2$, $F_4$, $E_8$.

The aim of this note is to present a simple argument which proves

{\bf Theorem 1.} {\em Let $G$ be an almost simple group. Then any cuspidal 
local system on $G$ is clean except, possibly, two cases: $p=2$ and $G$ is 
of type $F_4$ or $E_8$.}

Moreover, in each unsettled case there is exactly one cuspidal local system
which is not known to be clean. 
In view of results of Lusztig \cite{LuCS} and Shoji \cite{Sh1,Sh2} Theorem 1 
is new only for $G$ of type $E_6$ and $p=2$ and $G$ of type $E_7$ and $p=3$.

\subsection{Applications} Recall that Lusztig defined a class of perverse
sheaves on group $G$ called the admissible complexes, see \cite{LuIC}. It
is known that the character sheaves form a subset of admissible sheaves,
see \cite{LuCS}. It is known in many cases that actually character sheaves
coincide with admissible complexes, see \cite{LuCS, Sh1, Sh2}.

{\bf Theorem 2.} {\em The class of character sheaves coincides with the class 
of admissible complexes for any $p\ge 0$.}

In view of Lusztig's results in \cite{LuCS} 7.1 Theorem 2 is an immediate
consequence of Theorem \ref{main} (a) below.

Recall that Lusztig defined generalized Springer correspondence in \cite{LuIC}
which is a bijection between the irreducible $Ad(G)-$equivariant local systems
supported on the unipotent orbits and the irreducible representations
of some collection of Coxeter groups. The generalized Springer correspondence
is known explicitly \cite{LS, Spa, LuCS} in all cases with two very small 
gaps: in the case when the Coxeter group is of type $G_2$ there is an 
ambiguity in attaching the local systems to two-dimensional representations 
of this group for $G$ of type $E_6$ when $p=2$ and $E_8$ when $p=3$. This 
ambiguity can be now removed using the method used by Lusztig in \cite{LuCS} 
24.10 to handle a similar problem for $E_6$, $p>3$. We use below the 
notations from \cite{Spa}.

{\bf Proposition 1.} {\em Let $G$ be of type $E_6$ and $p=2$ (respectively,
of type $E_8$ and $p=3$). Under the generalized Springer correspondence
the reflection representation of $G_2$ corresponds to the local system
supported on the orbit of type $A_5$ (respectively $E_7$).}

We omit the proof since it coincides with \cite{LuCS} 24.10 (note that the
calculation of the corresponding generalized Green functions are almost 
identical for $E_6$ and $E_8$). 

\subsection{Acknowledgment} I learned the definition of the automorphism 
$\Theta_\cF$ which is crucial for this paper from Roman Bezrukavnikov who
in turn learned it from Vladimir Drinfeld. I am happy to thank both of them. 
I am deeply grateful to George Lusztig for very useful conversations. Thanks
are also due to Toshiaki Shoji for interesting comments. This work was 
supported in part by NSF grants DMS-0098830 and DMS-0111298.

\section{Proofs}
\subsection{} Let $a, p: G\times G$ be the adjoint action and the second
projection respectively: $a(g,x)=gxg^{-1}$, $p(g,x)=x$.
Let $\cF$ be a complex of (constructible) sheaves on $G$ which 
is $Ad(G)-$equivariant in the naive sense: we are given an isomorphism
$\xi : a^*\cF \to p^*\cF$ satisfying the cocycle relation, see e.g. \cite{Sh1}.
Let $\Delta: G\to G\times G$ be the diagonal embedding. Obviously
$p\Delta =a\Delta =\Id$. The following definition is crucial for this note:

\begin{definition} 
We define a canonical automorphism $\Theta_{\cF}$ of
$\cF$ as the composition:

$$\cF =(a\Delta)^*\cF =\Delta^*a^*\cF\stackrel{\Delta^*\xi}{\longrightarrow}
\Delta^*p^*\cF=(p\Delta)^*\cF=\cF.$$
\end{definition}

\begin{remark} (i) The definition above makes sense in the case when $G$ is
a finite group. In this case $\Theta_\cF$ is well known in the conformal
field theory under the name of $T-$matrix, see e.g. \cite{BK}.

(ii) The automorphism $\Theta_\cF$ was used in \cite{E,Sh1} to prove that
the characteristic functions of character sheaves are eigenvectors for 
Shintani descent.
\end{remark}

We are going to apply this definition for two kinds of complexes: the usual
constructible sheaves and the perverse sheaves.

\begin{lemma} \label{Schur}
Let $\cF$ be a simple $Ad(G)-$equivariant perverse sheaf. Then 
$\Theta_\cF =\theta_{\cF}\Id_\cF$ for some scalar $\theta_\cF$.
\end{lemma}

\begin{proof} This is an immediate consequence of the Schur's Lemma. 
\end{proof}

Let $\cO$ be an adjoint orbit in $G$ and let $x\in \cO$. Recall that the 
functor $\cL \mapsto \cL_x$ defines an equivalence
\{ $G-$equivariant local systems $\cL$ on $\cO$\} $\to$ \{ Representations
of $A_G(x):=Z_G(x)/Z_G(x)^0$\} . Let $\bar x$ be the class of $x\in Z_G(x)$ 
in the group $A_G(x)$. Observe that $\bar x\in A_G(x)$ is central.

\begin{lemma} \label{local}
Under the equivalence above we have $\Theta_\cL =\bar x|_{\cL_x}$.
\end{lemma}

\begin{proof} This is a direct consequence of definition. \end{proof}

Combining Lemmas \ref{Schur} and \ref{local} one can calculate the number
$\theta_\cF$ for an irreducible $Ad(G)-$equivariant perverse sheaf $\cF$
in the following way: take any point $x\in G$ such that $\cF_x\ne 0$, then
$\bar x\in A_G(x)$ acts on $\cF_x$ via the scalar $\theta_\cF$. 

Now let $P$ be a parabolic subgroup of $G$ with Levi quotient $L$.
Recall (see \cite{LuCS}) that the induction functors $\Ind_P:$ 
\{ $Ad(L)-$equivariant sheaves on $L$\} $\to$ \{ $Ad(G)-$equivariant sheaves 
on $G$\} is defined as follows: consider the variety $\tilde G_P=\{ (x,gP)\in
G\times G/P|g^{-1}xg\in P\}$. The group $G$ acts on $\tilde G_P$ in the
following way: $h\cdot (x, gP)=(hxh^{-1}, hgP)$. It is easy to see that we
have an equivalence: $\nu$: \{ $G-$equivariant sheaves on $\tilde G_P$\} 
$\simeq$ \{ $Ad(P)-$equivariant sheaves on $P$\} . Let $m: P\to L$ be the
canonical projection and let $n: \tilde G_P\to G$ be defined by $n(x,gP)=x$.
Then the functor $\Ind_P:=n_!\nu^{-1}m^*$ (we use here just naive notion of 
the equivariance but a similar construction holds for example in the
equivariant derived category).

\begin{lemma} \label{ind}
The automorphism $\Theta_\cF$ commutes with the induction functor:
$$\Ind_P(\Theta_\cF)=\Theta_{\Ind_P(\cF)}.$$\end{lemma} 

\begin{proof} Easy. \end{proof}

\begin{remark} A statement similar to Lemma \ref{ind} involving the Shintani
descent twisting operator instead of $\Theta_\cF$ is contained in \cite{E,Sh1}.
\end{remark}

\subsection{Calculation of $\theta_\cF$ for some unipotent cuspidal pairs}
We refer the reader to \cite{LuIC} for the definition of the cuspidal pair
$(C,\cL)$ for the group $G$ (recall only that here $C$ is some inverse image of
conjugacy class under projection $G\to G/Z_G^0$ and $\cL$ is some 
$Ad(G)-$equivariant local system on $G$). In this section we consider the
case when $G$ is semisimple and $C$ is an unipotent class (such
cuspidal pairs are called unipotent) and calculate $\theta_\cF$ for $\cF =\cL$
(extended by zero to $G$).  First note that

\begin{lemma}
Assume that the characteristic of $k$ is good for $G$. Then $\theta_\cL=1$
for any unipotent cuspidal pair $(C,\cL)$.
\end{lemma}

\begin{proof} Obviously it is enough to prove the Lemma for simply connected
almost simple groups.
For groups of type $A$ the order of $A_G(u)$ for any unipotent 
element $u$ is relatively prime to the characteristic of $k$ and hence 
$\bar u=1$ (since the order of $u$ is some power of the characteristic); the 
result follows from Lemma \ref{local}. For other classical groups $A_G(u)$ is 
a 2-group; thus $\bar u\in A_G(u)$ is trivial. Similarly, for groups of type
$E_6, E_7$ the order of $A_G(u)$ is always relatively prime with the
characteristic (assumed to be good). If the group $G$ is of type $G_2, F_4,
E_8$ then the unipotent cuspidal local system $(C\,cL)$ is unique. In these 
cases for $u\in C$ one has $A_G(u)=S_3,S_4,S_5$ respectively. Thus $A_G(u)$
has trivial center, hence $\bar u=1$.
\end{proof}

Now we are going to calculate $\theta_\cL$ for unipotent cuspidal pairs in
the exceptional groups. We are going to use the following
fact due to T.~Springer and B.~Lou \cite{Spr, Lou}:

\begin{theorem} \label{blou}
Let $u\in G$ be a regular unipotent element and let 
$U\subset G$ be the maximal unipotent subgroup containing $u$. Then $Z_G(u)=
Z_G\times Z_U(u)$; moreover the group $Z_U(u)/Z_U^0(u)$ is cyclic and is
generated by $\bar u$.
\end{theorem}

The list of unipotent cuspidal pairs for exceptional groups is given in 
\cite{Spa}. We give the values of $\theta_\cL$ in all bad characteristic
cases using the notations of {\em loc. cit.} In the table below $\zeta$ 
(respectively $\xi$) is a fixed primitive root of unity of degree 3 
(respectively 5). 

$$\begin{array}{|l|l|l|l|l|l|l|}\hline 
&G&char(k)&\mbox{class of} u&A_G(u)&\phi&\theta_\cL\\ \hline
1&G_2&2&G_2&\BZ_2&-1&-1\\
2&&&G_2(a_1)&S_3&\eps&1\\ \hline
3&G_2&3&G_2&\BZ_3&\zeta^{\pm 1}&\zeta^{\pm 1}\\
4&&&G_2(a_1)&\BZ_2&-1&1\\ \hline
5&F_4&2&F_4&\BZ_4&\pm i&\pm i\\
6&&&F_4(a_1)&\BZ_2&-1&-1\\
7&&&F_4(a_2)&D_8&\eps&1\\
8&&&F_4(a_3)&S_3&\eps&1\\ \hline
9&F_4&3&F_4&\BZ_3&\zeta^{\pm 1}&\zeta^{\pm 1}\\
10&&&F_4(a_3)&S_4&\eps&1\\ \hline
11&E_6&2&E_6&\BZ_6&-\zeta^{\pm 1}&-1\\
12&&&A_5+A_1&\BZ_6&-\zeta^{\pm 1}& 1\\ \hline
13&E_6&3&E_6&\BZ_3&\zeta^{\pm 1}&\zeta^{\pm 1}\\ \hline
14&E_7&2&E_7&\BZ_4&\pm i&\pm i\\ \hline
15&E_7&3&E_7&\BZ_6&-\zeta^{\pm 1}&\zeta^{\pm 1}\\
16&&&D_6(a_2)+A_1&S_3\times \BZ_2&-\eps&-1\\ \hline
17&E_8&2&E_8(a_1)&\BZ_4&\pm i&\pm i\\
18&&&E_7(a_2)+A_1&S_3\times \BZ_2&-\eps&-1\\
19&&&D_8(a_1)&D_8&\eps&1\\
20&&&2A_4&S_5&\eps&1\\ \hline
21&E_8&3&E_7+A_1&\BZ_6&-\zeta^{\pm 1}&\zeta^{\pm 1}\\
22&&&2A_4&S_5&\eps&1\\ \hline
23&E_8&5&E_8&\BZ_5&\xi^{\pm 1}, \xi^{\pm 2}&\xi^{\pm 1}, \xi^{\pm 2}\\ 
24&&&2A_4&S_5&\eps&1\\ \hline

\end{array}$$

\bigskip

\noindent
{\bf Comments on the calculation.} In cases 1,3,5,9,11,13,14,15,23 the 
calculation is immediate from \ref{blou}; in cases 2,4,7,8,10,19,20,22,24 the 
calculation is immediate from the fact that $\bar u$ is central in $A_G(u)$. 
The calculation in case 12 is as follows: assume that $\bar u\ne 1$ in this
case, then for the representation $\phi' =-1$ of $A_G(u)$ we will have 
$\phi'(\bar u)=-1$; this is a contradiction since the local system $\cL'$
corresponding to $\phi'$ appears in the principal series (see \cite{Spa})
and thus we have $\theta_{\cL'}=1$ by Lemma \ref{ind}; thus $\bar u=1$ and
the result follows. The similar method applies to cases 16,17,18,21. Finally,
for the case 6 see \cite{Sh1} 7.2 (it is stated there that $\bar u$ is
nontrivial in this case; one way to see this is an explicit calculation;
one can also use the results on Shintani descent, see {\em loc. cit.}).

\begin{remark} For types $G_2$, $F_4$, $E_8$ the numbers in the table were
computed by T.~Shoji \cite{Sh1, Sh2} as the eigenvalues of the twisting
operator.
\end{remark}

The important consequence of the calculation above is the following

\begin{corollary} \label{distinct}
Let $\cL_1\ne \cL_2$ be two cuspidal local systems on a 
simple group $G$. Then $\theta_{\cL_1}\ne \theta_{\cL_2}$ except for
$char(k)=2$ and $G$ is of type $F_4$ or $E_8$.
\end{corollary}

\begin{remark} Note that the results above together with the explicit knowledge
of the generalized Springer correspondence for exceptional groups \cite{Spa}
allow to determine $\bar u$ for all unipotent elements $u$ in these groups.
\end{remark}

\subsection{Cleanness of cuspidal local systems} Recall that a local system
$\cL$ on a locally closed subset $U$ of a variety $X$ is called clean if
$j_!\cL=j_{!*}\cL$ where $j: U\subset X$ is the obvious embedding.
It is expected that all cuspidal local systems are clean. Here is a main 
result of this note:

\begin{theorem} \label{main}
(a) Any cuspidal sheaf is a character sheaf.

(b) Assume that $G$ is an almost simple exceptional group. Let
$(C,\cL)$ be a cuspidal pair for $G$. Then $\cL$ is clean except possibly
two cases: $k$ is of characteristic 2 and

1) $G$ is of type $F_4$ and $C$ is unipotent orbit of type $F_4(a_2)$ (there
is a unique such cuspidal pair);

2) $G$ is of type $E_8$ and $C$ is unipotent orbit of type $D_8(a_1)$ (there
is a unique such cuspidal pair). 
\end{theorem}

\begin{proof} The theorem is known to be true for classical groups and for
exceptional groups in good characteristic \cite{LuCS}. The proof for 
exceptional groups in bad characteristic is quite similar to proofs in 
\cite{LuCS}. First one shows that any cuspidal character sheaf is clean 
(except, possibly, two cases in (b)) using Proposition 7.9 of \cite{LuCS} III 
with the action of the center replaced by the automorphism $\Theta_\cF$ and 
using Corollary \ref{distinct} (recall that it is enough to consider the 
cuspidal local systems supported on the unipotent orbits, see \cite{LuCS} 
7.11). Then the results of \cite{LuCS} provide a classification of 
character sheaves. Finally one compares the list of cuspidal character sheaves
with the list of cuspidal local systems (known from \cite{LuIC}) and deduces 
(a). The case when $p=2$ and $G$ is of type $F_4$ or $E_8$ requires 
additional arguments, see \cite{Sh1} 7.3 and \cite{Sh2} 5.3. 

\end{proof}

\end{document}